\documentclass[a4paper,reqno]{amsart}
\usepackage{amssymb}
\usepackage{latexsym}
\usepackage{amsmath}
\usepackage{euscript}
\usepackage{graphics,color}
\usepackage[all]{xy}
\usepackage[margin=3cm]{geometry}
\usepackage{MnSymbol} 

\newcommand{\be}{\begin{equation}}
\newcommand{\ee}{\end{equation}}
\newcommand{\ba}{\begin{eqnarray}}
\newcommand{\ea}{\end{eqnarray}}
\newcommand{\baa}{\begin{eqnarray*}}
\newcommand{\eaa}{\end{eqnarray*}}
\newcommand{\bb}{\begin{thebibliography}}
\newcommand{\eb}{\end{thebibliography}}

\newcommand{\bi}[1]{\bibitem{#1}}
\newcommand{\lab}[1]{\label{#1}}
\newcommand{\re}[1]{(\ref{#1})}

\newcounter{my}
\newcommand{\he}%
   {\stepcounter{equation}\setcounter{my}%
   {\value{equation}}\setcounter{equation}0%
   }%
\newcommand{\she}%
   {\setcounter{equation}{\value{my}}%
    }%

\renewcommand\t{\tilde}

\newcommand\vphi{\varphi}

\newcommand\ve{\varepsilon}

\newtheorem{pr}{Proposition}

\theoremstyle{definition}

\numberwithin{equation}{section}

\begin{document}

\title{The q-Heun operator of big q-Jacobi type and the q-Heun algebra}

\author{Pascal Baseilhac}
\author{Luc Vinet}
\author{Alexei Zhedanov}

\address{Institut Denis-Poisson CNRS/UMR 7013 - Universit\'e de Tours - Universit\'e d'Orl\'eans Parc de Grammont,
37200 Tours, FRANCE}
\address{Centre de recherches \\ math\'ematiques,
Universit\'e de Montr\'eal, P.O. Box 6128, Centre-ville Station,
Montr\'eal (Qu\'ebec), H3C 3J7}

\address{Department of Mathematics, School of Information, Renmin University of China, Beijing 100872,CHINA}

\begin{abstract}
The q-Heun operator of the big q-Jacobi type on the exponential  grid is defined. 	This operator is the most general second order q-difference operator that maps  polynomials of degree $n$ to polynomials of degree $n+1$. It is tridiagonal in bases made out of either q-Pochhammer or big q-Jacobi polynomials and is bilinear in the operators of the q-Hahn algebra. The extension of this algebra that includes the q-Heun operator as generator is described. Biorthogonal Pastro polynomials are shown to satisfy a generalized eigenvalue problem or equivalently to be in the kernel of a special linear pencil made out of two q-Heun operators. The special case of the q-Heun operator associated to the little q-Jacobi polynomials is also treated.
\end{abstract}



\maketitle

\section{Introduction}
\setcounter{equation}{0}
The purpose of this  paper is to introduce and characterize a q-difference 
analog of the Heun operator which we will call the q-Heun operator of the q-Jacobi type. In its general form,
 this operator will be shown to be related to the big q-Jacobi polynomials; the special case connected to the little q-Jacobi polynomials will also be distinguished.

This paper is the second of a series of articles dedicated to the study of the Heun-like equations that can be associated to the polynomials of the Askey scheme. The roots of this undertaking lie in the notion of algebraic Heun operator introduced in \cite{GVZ_band}. This operator which can be constructed in any bispectral context basically consists in the generic bilinear combination of the two operators providing the dual eigenvalue equations. Since all polynomials of the Askey scheme are bispectral, an algebraic Heun operator can be built for each family of the tableaux. The general name for these operators is motivated by the fact that the approach yields the standard Heun equation when applied to the Jacobi polynomials.

The difference analog of the standard Heun operator on the uniform linear lattice was obtained in \cite{VZ_HH}. We shall follow here a similar approach to deal with the exponential grid  and to obtain the q-analog of the Heun-Hahn operator presented therein. The main idea is first to find the most general operator defined on the given lattice that replicates on that space the property that the standard Heun operator has of mapping the polynomials of degree $n$ into polynomials of degree $n+1$. Second, one determines the corresponding algebraic Heun operator to find that the two approaches essentially give the same result. The bilinear construct is then easily seen to imply that the Heun operator in question acts tridiagonally on natural polynomial bases.

Examples of $q$-analogs of the Heun operators have recently been considered by Takemura in \cite{Takemura} where he introduces several types of such operators from degenerations of Ruijsenaars-Van Diejen operators \cite{Diejen}, \cite{Ruijsenaars}. In fact, one of those operators had already been introduced by Hahn \cite{Hahn_Heun}. The analysis that we offer in the present paper relates to two of these operators and provides an enlightening framework to understand their properties.

The definition of the algebraic Heun operator is intimately connected to a procedure
referred to as tridiagonalization \cite{IK1}, \cite{IK2}, \cite{GIVZ} which has been used to develop a bottom-up theoretical
description of orthogonal polynomials where families with higher number of parameters are
constructed and characterized from those with a lesser number of parameters. It is well known that the bispectral properties of the polynomials of the Askey
scheme can be encoded in algebras of quadratic type that are referred to as the Racah and Askey-Wilson
algebras in the case of the polynomials at the top of the hierarchies \cite{Zh_AW}, \cite{GLZ_Annals}. 
The algebras associated to the other families are similarly identified by the names of
the corresponding polynomials. Algebraically, the effect of tridiagonalization is hence to embed the more complicated algebras
into the simpler ones. 

It should be stressed that the most general tridiagonalization that leads to the algebraic Heun operator 
transcends of course the framework of orthogonal polynomials. In a recent study \cite{BMVZ}, we have shown for instance that the little q-Jacobi operator and a special tridiagonalization of it realizes the equitable presentation of the Askey-Wilson algebra. We shall supplement this by identifying the extension that results when the Heun operator is taken instead as one of the generators.

 It will also be observed that certain biorthogonal polynomials arise as solutions of the generalized eigenvalue problems associated to Heun-q-Jacobi difference operator.

The outline of the paper is as follows. In section 2 we construct the Heun operators on the exponential grid as second order q-difference operators that map polynomials of degree $n$ into polynomials of degree $n+1$. The operator stemming from the most general situation will be referred to as the big q-Heun operator. A special case that yields what we will call the little q-Heun operator will also be discussed. These two operators will be seen to coincide (up to trivial similarity and affine transformations) with two of the q-Heun operators introduced in \cite{Takemura}.

The q-Hahn algebra encoding the bispectrality of the q-Hahn polynomials will be exhibited in Section 3 as a special case of the Askey-Wilson algebra and seen to be realized by the q-difference and recurrence operators of the polynomials it is bearing the name.

The construction of the algebraic Heun operator corresponding the bispectral situation depicted by the q-Hahn algebra is carried out in Section 4 and the operator obtained in this way, will be identified with the big q-Heun operator. In the special case where one parameter is set to zero and tridiagonalization proceeds with the little q-Jacobi operator, the algebraic approach will yield expectedly the little q-Heun operator. With these observations made, it will be shown that the big q-Heun operator acts tridiagonally on big q-Jacobi  and q-Pochhammer polynomials.

The q-Heun algebra that results when the big q-Heun operator is combined with the q-Hahn difference operator is the focus of Section 5. Upon examining the reduction of this algebra to the q-Racah one, an Askey-Wilson triple will be naturally obtained.

It will be shown in Section 6 that the Pastro biorthogonal polynomials are in the kernel of a linear pencil made out of special q-Heun operators. Concluding remarks, in Section 7, will end the paper.

\section{Construction of the big q-Heun operator on the exponential grid}
\setcounter{equation}{0}

We will deal with second-order q-difference operators of the type
\be
W= A_1(x) T^+ + A_2(x) T^- + A_0(x) \mathcal{I}, \lab{gen_W_def} \ee
where 
\be
T^+ f(x) = f(qx), \; T^- f(x) = f(x/q) \lab{TPM_def} \ee
are q-shift operators. In general, such operators  act on infinite-dimensional spaces. Finite-dimensional restrictions are however possible and in those cases one can say that the operator $W$ acts on a finite q-exponential grid.

\vspace{3mm}

To have a parallel with  the ordinary Heun operator and the Heun-Hahn operator introduced in  \cite{VZ_HH}, we would like to find the operator $W$ of form \re{gen_W_def} such that:

{\bf (i)} $W$ is the most general second order q-difference operator  which sends any polynomial of degree $n$  to a polynomial of degree $n+1$.

\vspace{3mm}

In what follows we shall show that property (i) is equivalent to :

{\bf (ii)} $W$ is the most general second order differential operator which is tridiagonal with respect to the q-Pochhammer basis $\varphi_n(x)=(q^{-x};q)_n, \: n=0,1,2,\dots, $,
where  $(x;q)_n=(1-x)(1-xq) \dots (1-xq^{n-1})$,
\be
W \varphi_n(x) = \{\vphi_{n+1}(x), \vphi_n(x), \vphi_{n-1}(x)\}. \lab{3-diag_Pochh} \ee

\vspace{3mm}

Let us start with the  construction of the q-difference Heun operator meeting property (i). Assume that $W$  is generically given by  \re{gen_W_def}. Our first task is to determine the coefficients $A_i(x), \: i=0,1,2$. To do this we can act with $W$ successively on the first three monomials $1,x,x^2$.

Operating  on the constant yields:
\be
W\{1\} = A_1(x) +A_2(x) + A_0(x). \lab{W_const} \ee 
According to property (i), we find from \re{W_const} that 
\be
A_0(x) = \pi_1(x) -A_1(x) - A_2(x), \lab{A0_12} \ee
where $\pi_1(x)$ is a polynomial of first degree. Applying the operator $W$ to the function $x$,  we obtain similarly that
\be
A_1(x)  qx + A_2(x) x/q + A_0(x) x = \pi_2(x) \lab{2c_A12} \ee
with some polynomial $\pi_2(x)$ of second degree. Taking into account \re{A0_12} we get
\be
xA_1(x)(q-1)+xA_2(x)(q^{-1}-1) = \tilde \pi_2(x)  \lab{W_x^2} \ee
with $\deg(\t \pi_2(x)) \le 2$. Applying now $W$ to the monomial $x^2$, we see that the property (i) imposes:
\be
x^2 A_1(q^2-1) + x^2 A_2(x) (q^{-2}-1) = \t \pi_3(x) \lab{W_x^3} \ee
with $\deg(\t \pi_3(x)) \le 3$.  It follows from these results that $A_0(x), A_1(x), A_2(x)$ have the following general expressions
\be
A_1(x) = \frac{p_3(x)}{x^2}, \; A_2(x) = \frac{q p_3(x) +xp_2(x)}{x^2}, \; A_0(x) = -A_1(x) - A_2(x) + p_1(x) \lab{A012} \ee
with $\deg(p_i(x)) \le i$.

It is then easily verified that property (i) holds for any $n=3,4,\dots$. Indeed, the application of the operator $W$ to the monomial $x^n$ yields
\be
W x^n = (q^n-1)  \left\{ (1 -q^{1-n})p_3(x) -q^{-n}x p_2(x) \right\} x^{n-2} + p_1(x) x^n. \lab{W_x_n} \ee
It is clear from this expression that the rhs of \re{W_x_n} is a polynomial of degree not exceeding $n+1$. Hence expressions \re{A012} defines the most general second-order operators $W$ with the property (i). We will refer the operator $W$ thus obtained as the big q-Heun operator (the explanation of this name will be given in Section 4).

There is an interesting special case that occurs when $p_3(0)=0$, i.e. when $p_3(x) = x r_2(x)$ with some polynomial $r_2(x)$ of second degree. In this case the coefficients $A_1(x), A_2(x)$ simplify to 
\be
A_1(x) = \frac{r_2(x)}{x}, \quad A_2(x) = \frac{q r_2(x) +p_2(x)}{x}, \quad A_0(x) = -A_1(x) - A_2(x) + p_1(x) \lab{A012_little} \ee
It is seen that 
\be
A_k(x) = \frac{s_k(x)}{x}, \; k=0,1,2 \lab{A_k_little} \ee
with $s_k$ some second degree polynomials. These polynomials are arbitrary apart from the only condition that
\be
s_0(x)+s_1(x)+s_2(x)=p_1(x). \lab{s_cond} \ee
In this special case, we shall call $W$  the little q-Heun operator.

\vspace{2mm} 
We can now relate our results with those presented by Takemura in \cite{Takemura} where he identified several q-Heun operators arising as degenerations of Ruijsenaars--van Diejen operators (see \cite{Diejen}, \cite{Ruijsenaars}). The third q-Heun operator presented in \cite{Takemura} is 
\begin{gather}
 A^{\langle 3 \rangle} = x^{-1} \prod_{n=1}^3 \big(x- q^{h_n+1/2} t_n\big) T^- + x^{-1}
 \prod_{n=1}^3 (x- q^{l_n -1/2} t_n ) T^+ -\big(q^{1/2} +q^{-1/2} \big) x^2\nonumber\\
 \hphantom{A^{\langle 3 \rangle} =}{} +\sum _{n=1}^3 \big( q^{h_n} + q^{l_n} \big)t_n x + q^{(l_1 +l_2 +l_3 +h_1 +h_2 +h_3)/2} \big( q^{\beta/2} + q^{-\beta/2} \big) t_1 t_2 t_3 x^{-1} , \label{eq:qthird}
 \end{gather}
and the fourth q-Heun operator $A^{\langle 4 \rangle}$ given in \cite{Takemura} reads
\begin{gather}
 A^{\langle 4 \rangle} = x^{-1} \big(x-q^{h_1 + 1/2} t_1\big) \big(x-q^{h_2 +1/2} t_2\big) T^- \nonumber\\
\hphantom{A^{\langle 4 \rangle} =}{}
 + q^{\alpha _1 +\alpha _2} x^{-1} \big(x - q^{l_1 -1/2} t_1\big) \big(x - q^{l_2 -1/2} t_2\big) T^+ \nonumber \\
\hphantom{A^{\langle 4 \rangle} =}{} -\big\{ \big(q^{\alpha _1} +q^{\alpha _2} \big) x + q^{(h_1 +h_2 + l_1 +l_2 +\alpha _1 +\alpha _2 )/2} \big( q^{\beta/2} + q^{-\beta/2} \big) t_1 t_2 x^{-1} \big\}. \label{eq:qH}
\end{gather}

In these formulas $h_i,l_i, t_i, \alpha_i,\beta$ are parameters. It is the operator $A^{\langle 4 \rangle}$ that is the object of special attention in \cite{Takemura} as a q-generalization of the Heun operator.

As a matter of fact, an operator equivalent to $A^{\langle 4 \rangle}$ was first introduced by Hahn \cite{Hahn_Heun}; this q-Heun operator of Hahn can be presented in the form
\be
\frac{s_1(x)}{x} T^+ + \frac{s_2(x)}{x} T^- + \frac{s_0(x)}{x} \mathcal{I}, \lab{Hahn-q-Heun} \ee
where $s_0(x),s_1(x),s_2(x)$ are arbitrary quadratic polynomials. It is easily seen that indeed the operator \re{Hahn-q-Heun} coincide with the operator $A^{\langle 4 \rangle}$.

On the other hand, our little q-Heun operator corresponding to the special case \re{A_k_little}, has the same expression \re{Hahn-q-Heun} but with the additional restriction \re{s_cond}. We can eliminate this restriction using the similarity transformation 
\be
\t W = x^{\gamma} W x^{-\gamma} \lab{tWW} \ee
with some constant $\gamma$. Indeed, condition \re{s_cond} means that the  sum of the leading coefficients in $x^2$ of the polynomials $s_0(x), s_1(x), s_2(x)$ vanishes. It is clear that this condition can always be enforced by an appropriate choice of the parameter $\gamma$ in \re{tWW}. This indicates that condition \re{s_cond} is not an essential restriction and that we have indeed arrived at the most general $q$-Heun operators defined in \cite{Hahn_Heun} and \cite{Takemura}. 

A similar transformation can be used to obtain the q-Heun operator $ A^{\langle 3 \rangle}$. In the latter case, one should allow for scaling of the argument 
\be
\t x \to \ve x \lab{scaling_x} \ee
and the following "affine" transformations of the Heun operator 
\be
W \to \theta(x) \left(W + (\alpha x +\beta)\mathcal{I} \right) \lab{affine_W} \ee
with an arbitrary function $\theta(x)$ and arbitrary constants $\alpha, \beta$.
Operators related by such transformations are deemed equivalent. Indeed, the eigenvalue problem for the Heun operators is 
\be
W \psi(x) =0 \lab{W_psi} \ee  
and it is clear that multiplication of $W$ by any function $\theta(x)$ does not modify the equation \re{W_psi}. The term $(\alpha x +\beta)\mathcal{I}$ will change only the linear function $p_1(x)$ in \re{A012} and is hence not essential. The identification of our little q-Heun operator with $A^{\langle 4 \rangle}$ can then be made under this natural equivalence relation. We thus have: 
\begin{pr}
Up to similarity, scaling and affine transformations \re{tWW}, \re{scaling_x}, \re{affine_W}, the little q-Heun operator is equivalent to the operator $A^{\langle 4 \rangle}$ while the big q-Heun operator is equivalent to the operator $A^{\langle 3 \rangle}$.
\end{pr}

\section{The q-Hahn algebra}
\setcounter{equation}{0}
The Askey-Wilson (AW) algebra \cite{Zh_AW}, \cite{GLZ_Annals} consists of 3 generators $K_1,K_2,K_3$ with the commutation relations
\ba
&&[K_1,K_2]=K_3, \nonumber \\
&&[K_2,K_3]= r K_2 K_1 K_2 + \xi_1 \{K_1,K_2\} + \xi_2 K_2^2 + \xi_3 K_2 +\xi_4 K_1 + \xi_5 \mathcal{I}, \nonumber \\
&&[K_3,K_1]= r K_1 K_2 K_1 +\xi_1 K_1^2 + \xi_2 \{K_1,K_2\}+ \xi_3 K_1 + \xi_6 K_2 + \xi_7 \mathcal{I}, \lab{AWA}
\ea
where $\mathcal{I}$ is the identity operator and $r, \xi_i$ are structure parameters of the algebra. The Askey-Wilson algebra naturally describes the eigenvalue problems of the Askey-Wilson (q-Racah in the finite-dimensional case) polynomials and arises in various physical and mathematical contexts  where these polynomials appear (see \cite{Zh_AW}, \cite{GLZ_Annals},  \cite{Ter} for more details).

The parameter $r$ in front of the cubic terms $K_2K_1K_2$ and $K_1K_2K_1$ accounts for the q-deformation of the classical orthogonal polynomials. The parameter $q$ is defined as
\be
r=2-q-q^{-1} \lab{r_q} \ee
If $q=1$, the cubic terms disappear in \re{AWA} and we obtain the Racah algebra \cite{GLZ_Annals} which describes Racah (Wilson) polynomials.

The q-Hahn algebra \cite{Zh_AW} is a specialization of the Racah algebra when one of the parameters $\xi_1$ or $\xi_2$ becomes zero. Taking for example $\xi_2=0$, we obtain the relations 
\ba
&&[K_1,K_2]=K_3, \nonumber \\
&&[K_2,K_3]= r K_2 K_1 K_2 + \xi_1 \{K_1,K_2\}  + \xi_3 K_2 +\xi_4 K_1 + \xi_5 \mathcal{I}, \nonumber \\
&&[K_3,K_1]= r K_1 K_2 K_1 + \xi_1 K_1^2 + \xi_3 K_1 + \xi_6 K_2 + \xi_7 \mathcal{I}. \lab{HA}
\ea
The q-Hahn algebra describes the eigenvalue problems of the big q-Jacobi polynomials (q-Hahn polynomials for the finite-dimensional case).

There is a "canonical" realization of the q-Hahn algebra in terms of q-difference operators. Indeed, let us introduce the operator $X$ which is the multiplication by $x$:
\[
X=K_1 =x \lab{X_x} 
\]
and the operator $Y$ which is  the q-difference operator associated to the big q-Jacobi polynomials \cite{KLS}
\be
Y=K_2 = B(x) T^+ + D(x) T^- -(B(x)+D(x)) \mathcal{I}, \lab{Y_def}
\ee
where $T^{\pm}$ have the same meaning as in \re{TPM_def} and where the coefficients are \cite{KLS}:
\be
B(x) =aqx^{-2} (x-1)(bx-c), \quad D(x) = x^{-2}(x-aq)(x-cq). \lab{BD_big} 
\ee
Here $a,b,c$ are the parameters of the big q-Jacobi polynomials $P_n(x;a,b,c)$. The latter are eigenfunctions of the  operator $Y$ \cite{KLS}
\be
Y P_n(x) = \lambda_n P_n(x), \lab{Y_P}
\ee
with eigenvalues
\be
\lambda_n =(q^{-n}-1)(1-ab q^{n+1}). \lab{eig_H}
\ee
It is directly verified that the operators $X$ and $Y$ (together with their commutator $K_3=[X,Y]$) satisfy the q-Hahn algebra \re{HA} with
\ba
&&\xi_1=-(q-1)^2(ab+q^{-1}), \: \xi_3=(q-1)^2(ab+ac+a+c), \: \xi_4=(q-1)^2(ab-1)(q^{-1}-qab),  \nonumber \\
&& \xi_5=(q-1)^2(ab-1)\left(aq(b+c)-a-c \right), \: \xi_6=0, \; \xi_7= -ac(q-1)^2(q+1). \lab{xi_AW} 
\ea

The q-Hahn algebra encompasses a duality property. This means that in the basis of the big q-Jacobi polynomials, the operator $Y$ becomes diagonal \re{Y_P} while the operator $X$ becomes tridiagonal:
\be
X P_n(x) = P_{n+1}(x) + b_n P_n + u_n P_{n-1}(x), \lab{X_P_rec} 
\ee
where the recurrence coefficients $b_n,u_n$ are given in \cite{KLS}. We will not use their explicit expression in this paper.

\section{The algebraic q-Heun operator of the big q-Jacobi type}
\setcounter{equation}{0}
In \cite{GVZ_band} another approach to construct the Heun operator was proposed.  This is based on the notion of {\it algebraic q-Heun operator} defined  as the generic bilinear combination of the generators $X,Y$ of the algebra associated to the relevant bispectral problem, namely
\be
W = \tau_1 XY + \tau_2 YX + \tau_3 X + \tau_4 Y + \tau_0 \mathcal{I}. \lab{W_Heun}
\ee
(We are using generically the same symbol $W$ as before for the q-Heun operator. We shall in fact see that they coincide.)
Note that operators of such type were considered by Nomura and Terwilliger \cite{NT} in the finite-dimensional context of the Leonard pairs as the most general operators which are tridiagonal with respect to dual bases (that diagonalize either operator $X$ or $Y$).

The operators $W$ are called algebraic Heun operators because in the special case of the Jacobi algebra (where $Y$ is the Gauss differential hypergeometric operator) the operator $W$ coincides with the ordinary Heun operator (see \cite{GVZ_Heun} for details). In \cite{VZ_HH} it was shown that the Heun-Hahn operator is equivalent to the algebraic Heun operator with $X,Y$ belonging to the  Hahn quadratic algebra. It is therefore natural to expect that the big (and little) q-Heun operators will coincide with the algebraic Heun operator \re{W_Heun} when $X,Y$ are the generators of the $q$-Hahn algebra. 

Using the expression \re{Y_def} for the operator $Y$, we have for the algebraic Heun operator of the big q-Jacobi type
\be
W= A_1(x) T^+ + A_2 T^- + A_0(x) \mathcal{I} 
\lab{W_Heun_BDC} \ee
where
\be
A_1(x) = \frac{aq(x-1)(bx-c) \left( (\tau_1+q\tau_2)x + \tau_4\right)}{x^2},
\lab{A1_alg} \ee
\be
A_2(x) = \frac{(x-aq)(x-cq) \left( (\tau_1+q^{-1} \tau_2)x + \tau_4\right)}{ x^2},
\lab{A2_alg} \ee
\be
A_0(x) = -A_1(x)-A_2(x) + p_1(x),
\lab{A0_alg} \ee
with 
\ba
&&p_1(x) =\left( \tau_2 (q-1)(qab-q^{-1}) + \tau_3 \right) x +\tau_0 +  (1-q) \left(qa(b+c) -a-c \right)\tau_2.
\lab{p1_alg} \ea

Comparison of \re{A1_alg}-\re{A0_alg} with \re{A012} shows that these expressions have the same structure. Using the scaling transformation $x \to \delta x$ of the argument, it is possible to put one of the roots of the polynomial $p_3(x)$ equal to 1. It can be checked that the coefficients $A_0(x),A_1(x),A_2(x)$ in \re{W_Heun_BDC} are the same as in \re{gen_W_def}.  This means that the polynomials $p_3(x),p_2(x),p_1(x)$ in \re{A012} can be expressed in terms of the parameters $\tau_i, \: i=0,1,\dots,4$ and $a,b,c$ and vice versa. We thus established the following new characterization property:

\vspace{3mm}

{\bf (iii)} the big q-Heun operator \re{W_Heun_BDC} is the second-order q-difference operator  given by the bilinear Ansatz \re{W_Heun} for $X$ and $Y$ the generators of the q-Hahn (or big q-Jacobi) algebra.

\vspace{3mm}

If $c=0$, the coefficients $A_i(x)$ become the same as \re{A012_little}. In this case the operator $Y$ is the difference operator for the little q-Jacobi polynomials. This hence leads to the little q-Heun operator.

The definition  \re{W_Heun} has advantages for the algebraic analysis of the q-Heun operators because it separates the parameters of the operator $W$ into two sets: {\it internal}  parameters $a,b,c$ coming from the big q-Jacobi polynomials and the {\it external} parameters $\tau_i, \: i=0,1,\dots,4$ originating  from the bilinear combination \re{W_Heun}.

It is now obvious that the operator $W$ is tridiagonal in the basis of the big q-Jacobi polynomials:
\be
W P_n(x) = \xi_{n+1} P_{n+1}(x) +\eta_n P_n(x) +\zeta_n u_n P_{n-1}(x).
\lab{W_3_big} \ee 
The coefficients $\xi_n, \eta_n, \zeta_n$ can easily be calculated using \re{X_P_rec} and \re{W_Heun}
\ba
&&\xi_n = \tau_1 \lambda_{n-1} + \tau_2 \lambda_n + \tau_3, \: \zeta_n=\tau_2 \lambda_{n-1} + \tau_1 \lambda_n + \tau_3, \nonumber \\
&& \eta_n = (\tau_1+\tau_2) \lambda_n b_n + \tau_3 b_n + \tau_4 \lambda_n + \tau_0. \lab{rec_cf_W} 
\ea
The inverse property can also  be checked. Namely, if $W$ is the most general operator which is 3-diagonal with respect to the big q-Jacobi polynomials, then taking $n=0,1,2$ we arrive at coefficients $A_1(x),A_2(x), A_0(x)$ having expressions \re{A012}. This means that the big q-Heun operator has the additional characterization property:

\vspace{3mm}

{\bf (iv)} The operator $W$ is the most general second order difference operator on the exponential grid which is tridiagonal with respect to the big q-Jacobi polynomials $P_n(x;a,b,c)$:
\be
W P_n(x) = \{ P_{n+1}(x), P_n(x), P_{n-1}(x) \}. \lab{Hahn_W_3diag-1} \ee

There is also the equivalent property which was already mentioned in Section 2:

(ii) The big q-Heun operator is 3-diagonal in the  q-Pochhammer basis:
\be
W \vphi_n(x) = \{\vphi_{n+1}(x), \vphi_n(x), \vphi_{n-1}(x) \} \lab{3_Pochh} \ee
where
\be
\vphi_n(x) = (x;q)_n = (1-x)(1-qx) \dots (1-q^{n-1}x). \lab{vphi_def} \ee
Property (ii) follows directly from formula \re{W_Heun}. Indeed,  the operator $X$ acts on the basis $\vphi_n(x)$ as follows:
\be
X \phi_n(x) = q^{-n} \left(\vphi_n(x) - \vphi_{n+1}(x).  \right) \lab{X_phi} \ee
The operator $Y$ has also simple action in this basis:
\be
Y \phi_n(x) = \lambda_n \phi_n(x) + \mu_n \phi_{n-1}(x) \lab{Y_phi} \ee
with $\lambda_n$ as in \re{eig_H} and with $\mu_n$ expressed as
\be
\mu_n = (1-q^{-n})(aq^n-1)(cq^n-1). \lab{mu_n} \ee
Formulas \re{X_phi}, \re{Y_phi} and \re{W_Heun} imply obviously the property \re{3_Pochh}.
We thus see that the bilinear Ansatz \re{W_Heun} can explain many properties of the q-Heun operator.

Finally note that the big q-Jacobi polynomials admit a finite-dimensional restriction when $c=q^{-N-1}$ with $N$ a positive integer. In this case the argument $x$ takes the finite number of values \cite{KLS}
\be
x_s = q^{-s} \quad  s=0,1,\dots N \lab{x_qH} \ee
and the big q-Jacobi polynomials become  the q-Hahn polynomials \cite{KLS}. The big q-Heun operator \re{W_Heun_BDC} that results in this instance can then be called the Heun operator of q-Hahn type. This operator acts on the finite grid \re{x_qH}. In the limit $q \to 1$, the q-Hahn polynomials become the ordinary Hahn polynomials \cite{KLS} and the q-Heun operator turns into the Heun operator of Hahn type which acts on the uniform grid $x=0,1,\dots, N$. This operator was introduced and analyzed in \cite{VZ_HH}.

Summing up, we have established that the big q-Heun operator $W$ can be characterized by four equivalent properties (i)-(iv) as in the case of the Heun-Hahn operator  \cite{VZ_HH}.

\section{The Heun-AW algebra and its reduction to the q-Racah algebra}
\setcounter{equation}{0}
On the one hand, the pair of operators $X$ and $Y$ generates the q-Hahn algebra. On the other hand we have the new q-Heun operator $W$ which can  be constructed from the operators $X$ and $Y$ via the bilinear Ansatz \re{W_Heun} (tridiagonalization procedure). One can ask whether pair of the operators $W$ and $Y$ realize some algebraic relations similar to those of the Askey-Wilson (AW) algebra. Obviously, in general, the operators $Y$ and $W$ cannot simply yield the AW algebra, because the operator $W$ lies beyond the realm of "classical" operators with known explicit spectra. Nevertheless, it turns out that the operators $Y,W$ still satisfy relatively simple commutation relations which can be checked directly :
\be
[Y,[W,Y]] =    e_1 Y^3 + r YWY + s_1 Y^2 + s_2 \{Y,W\} + s_3 W + s_4 Y + s_5 \mathcal{I} \lab{RH_1}
\ee
and 
\be
[W,[Y,W]] = e_2 YWY + e_3 Y^3 + e_4 Y^2 + r WYW + s_6 W^2 + s_7 \{Y,W\} + s_8 W + s_9 Y + s_{10} \mathcal{I} . \lab{RH_2} 
\ee
In the relations above,   the terms with the coefficients $r$ and $s_i, \: i=1,\dots,7$ have structures similar to those of the AW algebra \re{AWA}. The  extra-terms with the coefficients $e_1, e_2,e_3,e_4$ are those that make  the algebra with relations \re{RH_1} and \re{RH_2} more general than the AW algebra. The same situation occurs with the ordinary Heun operator, where extra-terms appear with respect to the Racah  algebra  \cite{GVZ_Heun} and similarly also in the case of  the Heun-Hahn operator \cite{VZ_HH}. It is therefore natural to refer to the algebra \re{RH_1} - \re{RH_2}, as the Heun-AW algebra of the big q-Jacobi type. The coefficient $r$ has the same expression \re{r_q} as for the AW algebra. Expressions for the coefficients $s_i, \:i=1,2,\dots,10$ are rather complicated and we will not present them here; instead, we will concentrate on the extra-terms. The corresponding coefficients are
\ba
&&e_1=-\tau_4 r, \nonumber \\
&& e_2 = \tau_4 (1-q^{-1})(q^3-1), \nonumber \\ 
&&e_3 = -\tau_4^2 (1+q^2)(1-q^{-1})^2, \nonumber \\
&& e_4= \tau_4^2 (abq+1)(2-q+q^2) + \tau_4 (1+q+q^2) \left( \tau_0  +q(\tau_1+\tau_2)(a+c+ac+ab) \right) + \nonumber  \\
&&ac(q+1)(q^2+q+1) (q \tau_1+\tau_2)(\tau_1+q\tau_2). \lab{eee} \ea 
If all the extra terms vanish, the relations \re{RH_1}-\re{RH_2} become the ordinary relations for the AW-algebra. From \re{eee} it is visible that these terms vanish in 4 possible cases:

(i) $\tau_4=0, \; \tau_2=-q\tau_1$;

(ii) $\tau_4=0, \; \tau_2=-q^{-1} \tau_1$;

(iii) $\tau_4=0, \; c=0$;

(iv) $\tau_4=0, \; a=0.$

\noindent The case (iii) corresponds to tridiagonalization of the little q-Jacobi polynomials. This situation was considered in details in \cite{BMVZ}. It is an exact q-analog of the tridiagonalization of the ordinary hypergeometric operator discussed in \cite{GIVZ}.  We do not have particular observations for the case (iv). Let us however consider what happens with algebraic Heun operator $W$ in the cases (i) and (ii). For the case (i) we have
\be
W_1 = B_1(x) T^+ + x^{-1}v_1(x)\mathcal{I} \lab{i_W} \ee
with 
\be
B_1(x) = \frac{\tau_1 aq (1-q^2)(x-1)(bx-c)}{x}. \lab{B1} \ee
Similarly for the case (ii)
\be
W_2= B_2(x) T^- + x^{-1}v_2(x)\mathcal{I} \lab{ii_W} \ee
with
\be
B_2(x) = \frac{\tau_1 (1-q^{-2})(x-cq)(x-aq)}{x}. \lab{B2} \ee
In the above formulas $v_1(x)$ and $v_2(x)$ are some quadratic polynomials.

Irrespective of the association with the choices (i) and (ii), together the operators $W_1$, $W_2$ and $Y$ possess remarkable (q-) commutation relations. Indeed, with appropriate choices of the  parameters $\tau_0, \tau_1, \tau_3$ (they can be different for $W_1$ and for $W_2$) and by an affine transformation of the operator $Y \to \alpha Y + \beta$, one can reduce the commutation relations of these operators to those of the $\mathbb{Z}_3$ symmetric presentation of the Askey-Wilson algebra \cite{Ter}
\be
Y W_1 - qW_1 Y = W_2 + \omega_1, \; W_1W_2 - qW_2 W_1 = Y + \omega_2, \; W_2 Y - q Y W_2 = W_1 + \omega_3, \lab{Z_3_AW} \ee
thus showing that the operators $Y, W_1, W_2$ form an Askey-Wilson triple.

\section{Bispectrality of Pastro polynomials and q-Heun pencil}
\setcounter{equation}{0}
The Pastro polynomials were introduced in \cite{Pastro}. They have the following explicit expression
\be
P_n(z) = \kappa_n \: {_2}\Phi_1 \left( {q^{-n},b \atop \frac{b}{a} q^{1-n}} ; \frac{z b^2}{a}\right) \lab{Pastro} \ee
with two arbitrary parameters $a$ and $b$.
The coefficient $\kappa_n$ is needed to make them monic, i.e. to ensure that
\be
P_n(z) = z^n + O(z^{n-1}). \lab{P_monic} \ee
In contrast to the ordinary orthogonal polynomials, the Pastro polynomials satisfy the recurrence relation \cite{VZ_Pastro}
\be
P_{n+1}(x) + g_n P_n(x) = x \left(P_n(x) + e_n P_{n-1}(x)\right) \lab{rec_Pastro} \ee
with some coefficients $g_n, e_n$. Equation \re{rec_Pastro} can be considered as a generalized eigenvalue problem for two bi-diagonal matrices. Polynomials satisfying recurrence relations of the type \re{rec_Pastro} are called biorthogonal polynomials of Laurent type \cite{Zh_L}.  

The explicit expression \re{Pastro} of the Pastro polynomials is very similar to that of the little q-Jacobi polynomials \cite{KLS} (with a slightly different dependence on the number $n$). Hence one can derive the difference equation of the Pastro polynomials starting from the difference equation of the little q-Jacobi polynomials \cite{KLS}. Let us provide the results without giving the details of the calculations: the Pastro polynomials satisfy the generalized eigenvalue problem
\be
L_1 P_n(z) = \lambda_n L_2 P_n(z), \lab{Pastro_GEVP} \ee
where $L_1, \: L_2$ are q-difference operators
\be
L_1 = \left(\frac{b^2 x-q}{(q-1)x}\right) T^+ + \frac{q-bx}{(q-1)x}, \; L_2=\left(\frac{aq-b^2x}{bx}\right) T^- + \frac{b^3 x-aq}{bx}. \lab{L_12} \ee
The eigenvalue has the following simple expression
\be
\lambda_n= \frac{q^n}{q-1}. \lab{lambda_Pastro} \ee 
It is easy to see that both $L_1$ and $L_2$ are operators belonging to the little q-Heun family \re{i_W} and \re{ii_W}. Thus the q-Heun operator pencil 
\be
(L_1 - \lambda_n L_2 ) P_n(z) =0 \lab{pencil_Pastro} \ee
serves as the q-difference equation for the Pastro polynomials.  

\noindent In \cite{VZ_HH} it was demonstrated that some biorthogonal rational functions satisfy the linear pencil equation \re{pencil_Pastro} where $L_1$ and $L_2$ both belong to a Heun operator family arising from the Hahn algebra. Moreover, in \cite{MNR} the ordinary Heun operator is seen to play the crucial role in finding differential equations for some classical biorthogonal rational functions of the Chebyshev type.  We thus see that linear pencils of algebraic Heun operators form an important tool  in the construction of the difference equations for biorthogonal polynomials or rational functions.

\vspace{5mm}

\section{Conclusion}
\setcounter{equation}{0}
We have constructed a big q-Heun operator which possesses four properties ((i)-(iv)) analogous to those of the ordinary Heun and Heun-Hahn operators. 

It was found that the algebra generated by the operators $Y,W$ differs from the Racah algebra by four additional terms. This new algebra reduces to the AW- algebra for four possible specializations of the q-Heun operator which allow the identification of a q-Racah triple.

We have showed also that the linear pencil of two Hahn-Heun operators gives a solution of the difference equation for Pastro biorthogonal polynomials. 

It is expected that the generalization of the analysis presented here  to operators on the Askey-Wilson grid will produce the general q-Heun operator of the Askey-Wilson type. One can reasonably expect that the algebraic Heun operator associated to the AW-algebra will also possess properties analogous to the other characteristics that have been seen to prevail in the cases analyzed so far. One may also think that this will allow the identification with the other operators listed in \cite{Takemura}. We plan to report soon on this.

\bigskip\bigskip
{\Large\bf Acknowledgments}

PB and AZ would wish to ackowledge the hospitality of the CRM during the course of this work.  PB is supported by the CNRS.  The research
of L.V. is funded in part by a discovery grant from the Natural Sciences and Engineering Research
Council (NSERC) of Canada. Work of A.Z. is supported by the National
Science Foundation of China (Grant No.11771015).

\vspace{15mm}

\bb{99}


\bi{BMVZ} P.Baseilhac, X. Martin, L.Vinet and A.Zhedanov, {\it Little  and big q-Jacobi polynomials and  the Askey-Wilson algebra}, arXiv:1806.02656.

\bi{Diejen} van Diejen J.F., {\it Integrability of difference Calogero–Moser systems}, J. Math. Phys. {\bf 35} (1994), 2983-–3004.




\bi{GIVZ} V. Genest, M.E.H. Ismail, L. Vinet, A. Zhedanov {\it Tridiagonalization of the hypergeometric operator and the Racah-Wilson algebra}, arXiv:1506.07803.

\bi{GZ_preprint} Ya.I.Granovskii and A.Zhedanov, {\it Exactly solvable problems and their quadratic algebras}, Preprint, DonFTI, 1989.

\bi{GLZ_Annals} Ya. A. Granovskii, I.M. Lutzenko, and A. Zhedanov, {\it Mutual integrability, quadratic algebras,
and dynamical symmetry}. Ann. Phys. {\bf 217} (1992),  1--20.

\bi{GVZ_Heun} F.A.Gr\"unbaum, L.Vinet and A.Zhedanov, {Tridiagonalization and the Heun equation}, J.Math.Physics {\bf 58}, 031703 (2017), arXiv:1602.04840.

\bi{GVZ_band} F.A.Gr\"unbaum, L.Vinet and A.Zhedanov, {Algebraic Heun operator and band-time limiting},  arXiv:1711.07862.

\bi{Hahn_Heun} W. Hahn, {\it On linear geometric difference equations with accessory parameters}, Funkcial.
Ekvac. {\bf 14} (1971), 73–-78.

\bi{IK1} M. E. H. Ismail and E. Koelink. {\it Spectral analysis of certain Schr\"odinger operators}. SIGMA,
{\bf 8}: 61–-79, 2012.

\bi{IK2} M. E. H. Ismail and E. Koelink. {\it The J-matrix method}.
Adv.Appl.Math., {\bf 56}, 379--395, 2011.

\bibitem{KLS} R. Koekoek, P.A. Lesky, and R.F. Swarttouw. {\it Hypergeometric orthogonal polynomials and their q-analogues}. Springer, 1-st edition, 2010.

\bibitem{MNR} A.Magnus, F.Ndayiragije  and A.Ronveaux, {\it Heun differential equation satisfied by some classical
biorthogonal rational functions} (to be published); \\ https://perso.uclouvain.be/alphonse.magnus/num3/biorthclassCanterb2017.pdf

\bi{NT} K.Nomura and P.Terwilliger, {\it Linear transformations that are tridiagonal with respect to both eigenbases of a Leonard pair}, Lin.Alg.Appl.
{\bf 420} (2007), 198--207.  arXiv:math/0605316.

\bi{Pastro} P.I. Pastro, {\it Orthogonal polynomials and some q-beta integrals of Ramanujan}, J. Math. Anal. Appl. {\bf 112} (1985),
517--540.




\bi{Ruijsenaars} Ruijsenaars S.N.M., {\it Integrable $BC_N$ analytic difference operators: hidden parameter symmetries and eigenfunctions},
in New Trends in Integrability and Partial Solvability, NATO Sci. Ser. II Math. Phys. Chem.,
Vol. 132, Kluwer Acad. Publ., Dordrecht, 2004, 217–-261.

\bi{Takemura} K. Takemura {\it On $q$-deformations of Heun equation}, arxiv:1712.09564.

\bi{Ter} P.Terwilliger, {\it Two linear transformations each tridiagonal with respect to an eigenbasis of the other}, Lin.Alg.Appl. {\bf 330} (2001), 149–-203.




\bi{VZ_Pastro} L.Vinet and A.Zheanov, {\it Spectral Transformations of the Laurent Biorthogonal Polynomials, II. Pastro Poly-
nomials}, Canad. Math. Bull. Vol. 44 (3), (2001), 337--345.

\bi{VZ_HH} L.Vinet and A.Zhedanov, {\it The Heun operator of Hahn type},  arXiv:1808.00153.

\bi{Zh_AW} A. S. Zhedanov, {\it "Hidden symmetry" of Askey-Wilson polynomials}, Theoret. and Math. Phys. {\bf 89}
(1991), 1146–-1157.

\bi{Zh_L} A.Zhedanov, {\it The Classical Laurent biorthogonal polynomials}. J. Comp. Appl. Math. {\bf 98} (1998), 121--147.


\end{thebibliography}

\end{document}